\documentclass{elsart}

\usepackage{amsmath,amsbsy,amssymb}
\usepackage[usenames]{color}

\newcommand{\bb}[1]{{\mathbb #1}}
\newcommand{\ol}{\overline}
\newcommand{\lbl}[2]{}

\newcommand{\re}{{\rm Re}\ \!\!}

\newcommand{\bs}{\boldsymbol}

\newtheorem{theorem}{Theorem}[section]{}
{}
\newtheorem{proposition}{Proposition}[section]{}
{}

\theoremstyle{definition}
{}

\newenvironment{proof}[1][{\bf Proof.}]
{\par\noindent{\bf #1} }{{\ }\hfill{$\square$}}

\begin{document}

\begin{frontmatter}

\title{Rigidity of Analytic Functions at the Boundary}
\author{Mustafa Arslan}\\
\address{Louisiana State University, Department of Mathematics\\
Baton Rouge, LA 70803-4918}


\begin{abstract}
A new elementary proof for a theorem of D. Burns and S. Krantz on
the rigidity of the analytic self maps of the unit disc was recently
discovered by L. Baracco, D. Zaitsev, and G. Zampieri. We use their
argument to generalize Burns-Krantz theorems on the unit disc and on
the unit ball of $\bb C^n$.
\end{abstract}

\end{frontmatter}
\section{Introduction}
A theorem of D. Burns and S. Krantz (\cite[theorem 2.1]{BK}) states
that if an analytic self-map $f$ of the unit disc $D$ satisfies
\[f(z)=z+O(|z-1|^4)\ \ \ \ \mbox{on $D$}\] then $f(z)=z$. We will
refer to this theorem as Burns-Krantz theorem. The Burns-Krantz
theorem was generalized to finite Blaschke products by D. Chelst
\cite{C}. The proof of both Burns-Krantz theorem and the that of the
theorem of Chelst uses Hopf lemma. In a recent paper, \cite{BZZ}, L.
Baracco, D.Zaitsev and G. Zampieri gave a new and elementary proof
for the Burns-Krantz theorem (\cite[proposition 3.1.]{BZZ}). Using
this we are able to prove the following more general rigidity
theorem at the boundary.
\begin{theorem}\label{t1} Let
$\varphi$ be an analytic function on the unit disc $D$, having
exactly $n$ zeros, counting multiplicities, in $D$. Let $f$ be an
analytic function on $D$ such that $|f|\leq|\varphi|$ on $\partial
D$. Suppose there are $n$ distinct points $\tau_j$ on $\partial D$
and positive integers $m_j$, $j=1,\ldots,k$ with $\sum m_j=n+1$,
such that $\varphi$ is away from zero around $\tau_j$ and
\lbl{1}{.3}
\begin{equation}\label{1}
f(z)=\varphi(z)+o(|z-\tau_j|^{2m_j-1})\ \ \mbox{as $D\ni
z\rightarrow \tau_j$}
\end{equation}
for each $j$. Then $f\equiv\varphi$ on $D$
\end{theorem}\lbl{t1}{-.2}

The conditions of the theorem \ref{t1} already exist in the
Burns-Krantz theorem and Chelst's theorem as any finite Blaschke
product $f$ have modulus 1 on $\partial D$ and the size of the set
$f^{-1}(\tau)$, $\tau\in\partial D$, is equal to the number of zeros
of $f$ in $D$ (counting multiplicities).

Let $B^n$ denote the unit ball of $\bb C^n$ and $B_r^n(\bs \tau)$
the ball of radius $r>0$ around $\bs \tau\in\bb C^n$. Another
theorem of Burns and Krantz (\cite[theorem 3.1]{BK}) states that if
$\Phi$ is a holomorphic map from the unit ball $B^n$ of $\bb C^n$
($n\geq2$) into itself which satisfies
\[\Phi(\bs z)=\bs z+O(|\bs z-\bs 1|^3)\] on $B^n$ then $\Phi(\bs z)=\bs z$ for
all $\bs z\in B^n$.

A similar approach lets us generalize this in the following fashion.
\begin{theorem}\label{t2}
Let $\Phi$ be a holomorphic map from the unit ball $B^n$ of $\bb
C^n$ into $\bb C^m$ with polynomial components of total degree at
most $s$ and none of which vanish at $\bs 1$. If $F$ is another
holomorphic map on $B^n$ such that $|F(\bs z)|\leq|\Phi(\bs z)|$ for
all $\bs z\in\partial B^n\cap B^n_\varepsilon(\bs 1)$ for some
$\varepsilon>0$ and
\begin{equation}\label{2}
F(\bs z)=\Phi(\bs z)+o(|\bs z-\bs 1|^{2s+1}) \ \ \ \mbox{as
$B^n\ni\bs z\rightarrow \bs 1$}
\end{equation}\lbl{2}{-.4}
then $F\equiv \Phi$ on $B^n$.
\end{theorem}\lbl{t2}{-.2}
Note that in this theorem we require the strong condition
$|F|\leq|\Phi|$ to hold only on a neighborhood of $\bs 1$ in
$\partial B^n$.

\section{Proofs}
\subsection{Proof of Theorem \ref{t1}}
Let $\alpha_j$, $j=1,\ldots,n$ be all the roots of $\varphi$ in $D$
repeated, if necessary, to count the multiple roots. Then we can
write
\[\varphi=(z-\alpha_1)\cdots(z-\alpha_n)\cdot u(z)\]where $u(z)$ is
an analytic function on $D$ having no zeros on $D$. Put
$\psi=\varphi/u$ and $g=f/u$. Since $|g(z)|\leq|\psi(z)|$ and $\psi$
is a polynomial, $g$ is a bounded analytic function on $D$ and since
$u$ is away from zero around $\tau_j$,
\[\psi(z)-g(z)=o((z-\tau_j)^{2m_j-1})\ \ \ \mbox{as $D\ni z\rightarrow \tau_j$}\]
for each $j=1,\ldots,k$.  Let, for $\tau\in\bb C$ and $\varepsilon
>0$, $B_\varepsilon(\tau)$ be the disc with center $\tau$ and radius
$\varepsilon$. Set
\[B_\varepsilon=\bigcup_{j=1}^kB_\varepsilon(\tau_j)\]
Now, observe that
\begin{equation}\label{3}
\re\left(\ol{\psi(e^{i\theta})}\cdot\frac{\psi(e^{i\theta})-g(e^{i\theta})}
{\prod_j|e^{i\theta}-\tau_j|^{2m_j}}\right)\geq 0
\end{equation}\lbl{3}{-.6}
For all $\theta \in[0,2\pi]$ except for those values of $\theta$
corresponding to $\tau_j$. Now we will show that the integral of the
function in \eqref{3} on $\theta$ from $0$ to $2\pi$ is 0:
\begin{align}\label{4}
\int_{\!\!\!\!\!\!\!\!\!\!\!\!\!\tiny\begin{array}{c} \\ 0\leq\theta\leq2\pi\\
e^{i\theta}\not\in B_\varepsilon\end{array}}
\ol{\psi(e^{i\theta})}&\cdot\frac{\psi(e^{i\theta})-{g(e^{i\theta})}}
{\prod_j|e^{i\theta}-\tau_j|^{2m_j}}d\theta\nonumber\\
&=\int_{\!\!\!\!\!\!\!\!\!\!\!\!\!\tiny\begin{array}{c} \\ 0\leq\theta\leq 2\pi \\
e^{i\theta}\not\in B_\varepsilon\end{array}}
\left(\prod_j(-e^{i\theta}\tau_j)^{m_j}\right)\ol{\psi(e^{i\theta})}\cdot\frac{\psi(e^{i\theta})
-{g(e^{i\theta})}}
{\prod_j(e^{i\theta}-\tau_j)^{2m_j}}d\theta\nonumber\\
&=-i\left(\prod_j(-\tau_j)^{m_j}\right)\int_{\partial D\backslash
B_\varepsilon}z^n\ol{\psi(z)}\cdot\frac{\psi(z)-{g(z)}}
{\prod_j(z-\tau_j)^{2m_j}}dz.
\end{align}\lbl{4}{-.7}
{For $z\in\partial D$, $z^n\ol{\psi(z)}=(1-\ol\alpha_1
z)\cdots(1-\ol\alpha_nz)$ which is the boundary function of the
analytic function \[\psi_1(z)=\prod_{j=1}^n(1-\ol\alpha_jz)\ \ \
z\in D.\] Now, ignoring the constant, we can rewrite integral in
\eqref{4} as
\[\int_{\partial D\backslash
B_\varepsilon}\psi_1(z)\cdot\frac{\psi(z)-{g(z)}}
{\prod_j(z-\tau_j)^{2m_j}}dz=\int_{\partial B_\varepsilon\cap
D}\psi_1(z)\cdot\frac{\psi(z)-{g(z)}} {\prod_j(z-\tau_j)^{2m_j}}dz\]
because the integrand is analytic in $D$. The last integral
decomposes into the sum
\[\sum_{j=1}^k\int_{\partial
B_\varepsilon(\tau_j)\cap D}\psi_1(z)\cdot\frac{\psi(z)-{g(z)}}
{\prod_{j=1}^k(z-\tau_j)^{2m_j}}dz\] and because of the condition
\eqref{1} each integral in the sum goes to zero as
$\varepsilon\rightarrow 0$.} This means that $|\psi|=\re(\ol \psi
g)$ almost everywhere on $\partial D$ but since $|\psi|\geq|g|$ on
$\partial D$ this is possible only when $\psi=g$ almost everywhere
on $\partial D$. Since both $\psi$ and $g$ are bounded functions,
$\psi\equiv g$ on $D$. Therefore $\varphi\equiv f$ on $D$.

Now we will show that the boundary condition given by the equation
\eqref{1} is the best possible. To be precise, for a given bounded
analytic function $\varphi$ on $D$ which is away from zero near the
boundary $\partial D$ and having zeros $\alpha_j\in D$ with
multiplicities $s_j\in\bb Z_+$, $j=1,2,\ldots,k$, we will construct
a function $f$ analytic on $D$ satisfying
$|f(e^{i\theta})|<|\varphi(e^{i\theta})|$ for almost all
$\theta\in[0,2\pi]$ and for each $\tau\in\partial D$, if we put
$\nu_\tau$ for the sum of those $s_j$ for which
$\frac{1+\alpha_j}{1+\ol\alpha_j}=\tau$ then
\begin{equation}\label{5}
f(z)=\varphi(z)+O(|z-\tau|^{2\nu_\tau})\ \ \ \mbox{for $z\in D$}
\end{equation}\lbl{5}{-.5}
and
\[f(z)=\varphi(z)+O(|z-1|^{2\nu_1+1})\ \ \ \mbox{for }z\in D.\]

The biggest part of the problem is covered by the following
proposition:
\begin{proposition}
For $\alpha_j\in D-\{0\}$, $s\in\bb Z_{\geq 0}$ and $s_j\in\bb Z_+$,
$j=1,\ldots,k$ we have
\begin{align}\label{5.5}
\!\!\!\!\!\!\!\!\!\!\!\!\!\!\!\!\!\!\!\!\!\!\!\!\left|z^s\prod_{j=1}^k(z-\alpha_j)^{s_j}
+ h(z)
(z-1)^{2s+1}\prod_{j=1}^k\left(z-\frac{1+\alpha_j}{1+\ol\alpha_j}
\right)^{2s_j}\right|&\nonumber\\
<\left|z^s\prod_{j=1}^k(z-\alpha_j)^{s_j}\right|\!\!\!\!\!\!\!\!\!
\!\!\!\!\!\!\!\!\!&
\end{align}\lbl{5.5}{-.6}
for almost all $z\in\partial D$, where
$h(z)=(-1)^nc\prod_j\frac{(1+\ol\alpha_j)^{2s_j}}{(1-\ol\alpha_jz)^{s_j}}$,
$n=s+\sum_js_s$ and $0<c<1/2^{-2n}$.
\end{proposition}
\begin{proof}
For $\alpha\in D$ we denote by $\gamma_\alpha$ the automorphism
$z\mapsto\frac{z-\alpha}{1-\ol\alpha z}$. Put $n$ for $s+\sum_js_j$.
Then for $z\in\partial D$ and $c>0$
\begin{align}\label{6}
\left|z^s\prod_j\right.\gamma_{\alpha_j}^{s_j}&+(-1)^nc(z-1)^{2s+1}\left.
\prod_j(\gamma_j -1)^{2s_j} \right|\nonumber \\
=&1+2(-1)^nc\ \re\left((\ol z(z-1)^2)^s(z-1)
\prod_j(\ol\gamma_{\alpha_j}(\gamma_{\alpha_j}-1)^2)^{s_j} \right)\nonumber \\
&\ \ \ +c^2|z-1|^{4s+2}\prod_j|\gamma_{\alpha_j}-1|^{4s_j}\nonumber \\
=&1-2^{n+1}c(1-\re z)^{s+1}\prod_j(1-\re\gamma_{\alpha_j})^{s_j}\nonumber\\
&\ \ \ +2^{2n+1}c^2(1-\re
z)^{2s+1}\prod_j(1-\re\gamma_{\alpha_j})^{2s_j}\\
\end{align}\lbl{6}{-.6}
Observe that if $c<2^{-2n}$, \eqref{6} is less than 1. So we have
\begin{equation}\label{7}
\left|z^s\prod_j\gamma_{\alpha_j}^{s_j}+(-1)^nc(z-1)^{2s+1}
\prod_j(\gamma_j -1)^{2s_j}
\right|\leq\left|z^s\prod_j\gamma_{\alpha_j}^{s_j}\right|\ \ \
\forall z\in\partial D.
\end{equation}\lbl{7}{-.6}
Multiplying both sides of \eqref{7} by
$\prod_j(1-\ol\alpha_jz)^{s_j}$ we obtain \eqref{5.5}.
\end{proof}

Now write
$\varphi(z)=u(z)(z-\alpha_1)^{s_1}\cdots(z-\alpha_k)^{s_k}$ where
$u(z)$ is non-vanishing analytic on $D$ (note that, here, we allow
any one $\alpha_j$ to be 0). It is now easy to see that the
function,
\[f(z)=u(z)\left(\prod_{j=1}^k(z-\alpha_j)^{s_j} +
h(z)\cdot(z-1)\prod_{j=1}^k\left(z-\frac{1+\alpha_j}{1+\ol\alpha_j}
\right)^{2s_j}\right)\] meets the requirements, where $h$ is as in
the proposition.
\subsection{Proof of Theorem \ref{t2}}\label{s2.1}
In this section, $\langle\cdot,\cdot\rangle$ will denote the
standard inner product (conjugate linear in the first variable) on
$\bb C^k$ for any dimension $k$ and $|\cdot|$ will be the
corresponding norm.

In the proof, we use the argument of the proof of \cite[proposition
3.1]{BZZ} again.

For a complex $(n-1)$-tuple $\alpha=(\alpha_2,\ldots,\alpha_n)$, set
\[Z_\alpha=\{\bs z=(z_1,\ldots,z_n)
\in\bb C^n:z_j=\alpha_j(z_1-1),\ j=2,3,\ldots n\}.\]
For each $\alpha\in\bb C^{n-1}$ the set $Z_\alpha\cap B^n$ is a disc
which has the following parametric description,
\[Z_\alpha\cap
B^n=\left\{\left(\frac{\zeta+|\alpha|^2}{1+|\alpha|^2},\frac{\alpha_2\zeta
-\alpha_2}{1+|\alpha|^2},\cdots,\frac{\alpha_n\zeta
-\alpha_n}{1+|\alpha|^2}\right):\zeta\in D\right\}.\] So if
$|\alpha|>2/\varepsilon=R$ then $Z_\alpha\cap B^n\subset
B^n_\varepsilon(\bs 1)$ with $\partial(Z_\alpha\cap
B^n)\subset\partial B^n\cap B^n_\varepsilon(\bs 1)$.

Fix an $\alpha$ with $|\alpha|>R$. Define functions $\varphi$ and
$f$, for $\zeta \in D$, by
\begin{eqnarray*}
\varphi(\zeta)&=&\Phi\left(\frac{\zeta+|\alpha|^2}{1+|\alpha|^2},
\frac{\alpha_1\zeta-\alpha_1}{1+|\alpha|^2},\ldots,
\frac{\alpha_n\zeta-\alpha_n}{1+|\alpha|^2}\right)\ \ \mbox{and}\\
f(\zeta)&=&F\left(\frac{\zeta+|\alpha|^2}{1+|\alpha|^2},
\frac{\alpha_1\zeta-\alpha_1}{1+|\alpha|^2},\ldots,
\frac{\alpha_n\zeta-\alpha_n}{1+|\alpha|^2}\right).
\end{eqnarray*}
Clearly, both $\varphi$ and $f$ are holomorphic maps from the unit
disc into $\bb C^n$ satisfying $|f|\leq|\varphi|$ on $\partial D$
and \[f(\zeta)=\varphi(\zeta)+o(|\zeta-1|^{(2s+1)})\ \ \ \mbox{as
$D\ni\zeta\rightarrow 1$}.\]So we have
\begin{equation}\label{8}
\re\left\langle\varphi(e^{i\theta}),\frac{\varphi(e^{i\theta})-
f(e^{i\theta})}{|e^{i\theta}-1|^{2s+2}}\right\rangle\geq 0\ \ \
\mbox{for all $\theta\in[0,2\pi]$}
\end{equation}\lbl{8}{-.5}

We repeat the same argument as in \ref{s2.1} on
\[\int_{\!\!\!\!\!\!\!\!\!\!\!\!\!\tiny\begin{array}{c} \\ 0\leq\theta\leq 2\pi \\
e^{i\theta}\not\in
B_\varepsilon(1)\end{array}}\left\langle\varphi(e^{i\theta}),\frac{\varphi(e^{i\theta})-
f(e^{i\theta})}{|e^{i\theta}-1|^{2s+2}}\right\rangle d\theta\] to
show that it is approaches to zero as $\varepsilon\rightarrow 0^+$
and conclude that $f\equiv\varphi+\eta$ for some holomorphic curve
$\eta$ which is normal to $\varphi$ on $\partial D$. But this
contradicts with the fact that $|f|\leq|\varphi|$ on $\partial D$
unless $\eta\equiv 0$. So we must have $f\equiv\varphi$ on $D$ (or
equivalently $F\equiv \Phi$ on $Z_\alpha\cap B^n$).

Now it is enough to show that the set
\[\bigcup_{|\alpha|>R}(Z_\alpha\cap B^n)\] contains an open subset of
$B^n$, but this is evident as
\[\bigcup_{|\alpha|>R}(Z_\alpha\cap B^n)=\left(\bigcup_{|\alpha|>R}
Z_\alpha-\{\bs 1\}\right) \cap B^n\] which is the nonempty
intersection of two open subsets of $\bb C^n$. This concludes the
proof of the theorem \ref{t2}.

The argument of L. Baracco, D. Zaitsev and G. Zampieri proves to be
very fruitful from which one can deduce several types of rigidity
theorems. One of these is the following, which can easily be proven
using the above setting.
\begin{proposition}
Let $\varphi$ be a holomorphic function on $B^n$ which is away from
0 near the point $\bs 1$. Let $f$ be another holomorphic function on
$B^n$ such that $|f(\bs z)|\leq|\varphi(\bs z)|$ for $\bs z\in
\partial B^n$ near $\bs 1$ and
\[f(\bs z)=\varphi(\bs z)+o(|\bs z-\bs 1|)\ \ \ \mbox{as $B^n\ni \bs z\rightarrow\bs 1$}.\]
Then $f\equiv \varphi$.
\end{proposition}


\begin{thebibliography}{10}
\bibitem{BZZ} L. Baracco, D. Zaitsev, G. Zampieri, {\it ``A
Burns-Krantz Type Theorem For Domains With Corners"}
arxiv:math.CV/0505264.
\bibitem{BK} D. Burns, S. Krantz, {\it ``Rigidity of Holomorphic Mappings and a
New Schwarz Lemma at the Boundary"} J. of Amer. Math. Soc. {\bf 7} (1994) no.
3, 661-676.
\bibitem{C} D. Chelst, {\it ``A Generalized Schwarz Lemma at the Boundary"}
{\bf 129} (2001) no.11, 3275-3278.
\end{thebibliography}
\end{document}